\documentclass[12pt]{amsart}
\usepackage{amsmath,amssymb, euscript, graphicx, hyperref}

\newtheorem{thm}{Theorem}[section]
\newtheorem{lem}[thm]{Lemma}
\newtheorem{prop}[thm]{Proposition}
\theoremstyle{definition}
\newtheorem{defn}[thm]{Definition}

\newtheorem{rem}[thm]{Remark}
\newtheorem{cor}[thm]{Corollary}

\begin{document}

\address{Azer Akhmedov, Department of Mathematics,
North Dakota State University,
Fargo, ND, 58102, USA}
\email{azer.akhmedov@ndsu.edu}

 \address{Warren Shreve, Department of Mathematics,
North Dakota State University,
Fargo, ND, 58102, USA}
\email{warren.shreve@ndsu.edu}

\begin{center} {\bf BALANCE IN RANDOM TREES} \end{center} 

\vspace{1cm}

\begin{center} {\bf Azer Akhmedov, \ Warren Shreve} \end{center}

\vspace{0.7cm}

 {\Small ABSTRACT: We prove that a random labeled (unlabeled) tree is balanced. We also prove that random labeled and unlabeled trees are strongly $k$-balanced for any $k\geq 3$. 
 
 {\em Definition:} Color the vertices of graph $G$ with two colors.  Color an edge
    with the color of its endpoints if they are colored with the same
    color.  Edges with different colored endpoints are left
    uncolored.  $G$ is said to be {\it balanced} if neither the
    number of vertices nor and the number of edges of the two
    different colors differ by more than one.
}

\section {Introduction}

    The notion of a balanced graph is defined [LLT] as follows:
  
  \medskip
  
 \begin{defn} Let $G = (V,E)$ be a finite simple graph, $k\geq 2$ be an integer, $c:V\rightarrow \{1, \ldots , k\}$ be a map. For all $i\in \{1, \ldots , k\}$, we write $V_i(c) = c^{-1}(\{i\}), E_i(c) = \{uv\in E \ | \ u,v\in V_i(c)\}$. We also write $v_i(c) = |V_i(c)|, e_i(c) = |E_i(c)|$. The map $c$ is called {\em a coloring}.
 \end{defn}
 
 \medskip
 
  The case of $k = 2$ is especially interesting. In this case, the sets $V_1(c), V_2(c), E_1(c), E_2(c)$ are called the sets of black vertices, white vertices, black edges, and white edges respectively. If the coloring $c$ is fixed we may drop it in the notation.

  \bigskip
  
  \begin{defn} A finite simple graph $G = (V,E)$ is called {\em balanced} if there exists a coloring $c:V\rightarrow \{1,2\}$ such that $|v_1(c)-v_2(c)|\leq 1$ and $|e_1(c)-e_2(c)|\leq 1$. A map $c:V\rightarrow \{1,2\}$ satisfying this condition is called {\em a balanced coloring}.
  \end{defn}
  
  \bigskip
  
  The graph in Figure 1. is balanced since we have shown the balanced coloring of it.
   
   \begin{figure}[h!]
  \includegraphics[width=3in,height=3in]{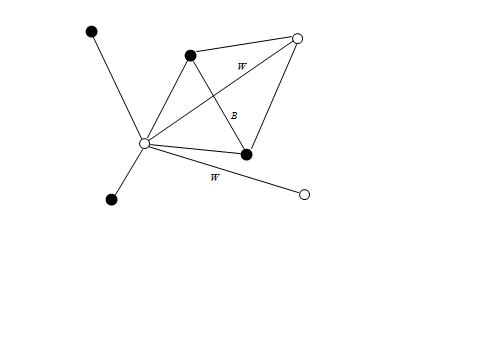}
\caption{with the given coloring, the graph has 4 black and 3 white vertices; it also has 2 white edges (labeled with a ``W") and 1 black edge (labeled with a ``B")}
\label{labelname}
\end{figure}

\bigskip

   It is not difficult to see that:
  
  \medskip
  
  a) The complete graph $K_n$ is balanced iff $n\leq 3$ or $n$ is even.
  
  \medskip
  
  b) The star $S_n$ is balanced iff $n\leq 5$; see Fig.2 for a balanced coloring of $S_5$.
  
  \medskip
  
  c) The double star $S_{p,q}$ is balanced iff $|p-q|\leq 3$.

  \medskip
  
  \begin{figure}[h!]
  \includegraphics[width=3in,height=3in]{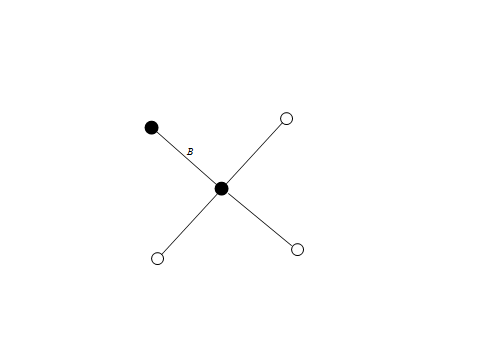}
\caption{a balanced coloring of $S_5$}
\label{labelname}
\end{figure}
  
  \medskip

  There are several other ways of defining balanced graphs which are similar to the one above but the definition we are using is the most interesting and challenging from the point of view of our problem.
   
   \medskip
  
	 In [Cah1], the author introduces a somewhat similar notion of a cordial graph, a generalization of both graceful and harmonious graphs. It has been conjectured by A.Rosa, G.Ringel and A.Kotzig that every tree is graceful ({\em Graceful Tree Conjecture}, [Ga]), and it has been conjectured by R.Graham and N.Sloane that every tree is harmonious (see [GS]). While these conjectures are still open, in [Cah2] it is proved that every tree is cordial.

	\medskip
	
  Not every tree is balanced; in this paper, we will be interested in the property of being balanced for a random labeled and unlabeled tree, as well as for random labeled graphs. 
  
  \medskip
  
  The main results of the paper are Theorem A and Theorem B stated below.
  
  \medskip
  
  {\bf Theorem A.} A random labeled (unlabeled) tree is balanced; more precisely, if $t_n (\tau_n)$ denotes the number of all labeled (unlabeled) trees on $n$ vertices,  and $b_n' (b_n'')$ denotes the number of all balanced labeled (unlabeled) trees on $n$ vertices, then $\displaystyle \lim _{n\rightarrow \infty }\frac{b_n'}{t_n} = 1$ and $\displaystyle \lim _{n\rightarrow \infty }\frac{b_n''}{\tau _n} = 1$. 
  
  \medskip
  
  \begin{rem} In this paper, for simplicity, we consider only uniform models of random graphs and random trees. The results can be extended to a large class of non-uniform models as well. Note that $t_n = n^{n-2}$(see [C] or [W]) and $\tau _n\sim C\alpha ^nn^{-5/2}$ for some positive constants $C$ and $\alpha $ (see [O]). 
  \end{rem}
  
  \medskip
  
  We also would like to introduce the notion of $k$-balanced graphs. 
  
 \medskip
 
  \begin{defn} Let $k\geq 2$. A finite simple graph $G = (V,E)$ is called $k$-{\em balanced} if there exists a coloring $c:V\rightarrow \{1, 2, \ldots , k\}$ such that $|v_i(c)-v_j(c)|\leq 1$ and $|e_i(c)-e_j(c)|\leq 1$ for all distinct $i,j\in \{1, 2, \ldots , k\}$. The map $c$ will be called a {\em $k$-balanced coloring}.
 \end{defn}
 
 \medskip
  
  \begin{defn} Let $k\geq 2$. A finite simple graph $G = (V,E)$ is called {\em strongly $k$-balanced} if there exists a coloring $c:V\rightarrow \{1, 2, \ldots , k\}$ such that $|e_i(c)| = 0, 1\leq i\leq k$, and $|v_i(c)-v_j(c)|\leq 1$ for all distinct $i,j\in \{1, 2, \ldots , k\}$. The map $c$ will be called a {\em strongly $k$-balanced coloring}.
 \end{defn}
 
 \medskip
  
  In more popular terms, a finite simple graph is strongly $k$-balanced iff it is {\em $k$-equitably colorable}.  In Section 5 we study some basic properties of $k$-balanced graphs. We prove the following theorem.
  
  \medskip
  
  {\bf Theorem B.} For all $k\geq 3$, a random (labeled) tree is strongly $k$-balanced.
  
 \medskip
 
  \begin{rem} Let us emphasize that Theorem B is orginally due to B.Bollob\'{a}s and R.Guy (see [BG]). Our proof in this paper is very different with some ingredients which might be interesting independently.
  \end{rem}
  
  \medskip
  
 \begin{rem}\label{rem:clique} It has been proved by I.Ben-Eliezer and M.Krivelevich (see [BK]) that a random graph is balanced. For $k\geq 3$, it seems quite plausible that a random graph is indeed $k$-balanced. However, notice that the clique number of a random graph on $n$ vertices is at least $2\log_2(n)$ (see [B]) thus a random graph is not strongly $k$-balanced.
  \end{rem}
  
  \bigskip
  
  {\em Acknowledgment:} We thank B.Gittenberger for the discussion and for bringing the reference [ES] to our attention. We are grateful to M.Krivelevich for bringing [BG] to our attention. We also would like to thank to I.Pak for his comments.      
  
  \medskip
  
  {\em Notes:} {\bf 1.} For any finite simple graph $G$, we will denote the maximal degree of $G$ by $d_{max}(G)$. 
  
  {\bf 2.} A vertex of degree one will be called a {\em leaf vertex} or simply a {\em leaf}. A non-leaf vertex $v$ is called a pre-leaf vertex if it is adjacent exactly to $m-1$ leaves where $m = deg(v)$. A pre-leaf vertex of degree two is called {\em special}.
  
  {\bf 3.} For $n\geq 2$, there exists a unique tree up to isomorphism with $n$ vertices and maximal degree at most two; we will call this tree a {\em string on $n$ vertices}. 
  
  {\bf 4.} For a tree  $G = (V, E)$ and a non-leaf vertex $v\in V$, a subset $A\subseteq V$ will be called a {\em branch} of $G$ with respect to $v$ if there exists a vertex $u$ adjacent to $v$ such that $A = \{x\in V \ | \ d(x,u) < d(x,v)\}$ where $d(.,.)$ denotes the distance in the tree $G$.
   
   \bigskip
   
  \section {Characterization of Balanced Graphs}
  
  In this section we observe some basic facts on balanced and $k$-balanced graphs. Let us first prove a very simple lemma which provides a necessary and sufficient condition for a graph to be balanced. 
  
  \medskip
  
  \begin{lem}\label{lem:trivial} Let $G$ be a finite simple graph with $n$ vertices, and degrees $d_1, \ldots , d_n$. $G$ is balanced if and only if there exists a partition $\{1,\ldots , n\} = I\sqcup J$ such that
  
  (i) $|$Card$(I)-$Card$(J)|\leq 1$
  
  (ii) $|\displaystyle \sum _{k\in I}d_k - \displaystyle \sum _{k\in J}d_k|\leq 2$
  \end{lem}
  
  \medskip
  
  {\bf Proof}. Let $G = (V,E), V = \{v_1, \ldots, v_n\}, \deg(v_i) = d_i, 1\leq i\leq n$.
  
  \medskip
  
  Assume $G$ is balanced with a balanced coloring $c:V\rightarrow \{1,2\}$. 
  
  \medskip
  
  Let $I = \{i \ | \ 1\leq i\leq n, c(v_i) = 0\}, J =  \{i \ | \ 1\leq j\leq n,   c(v_i) = 1\}$.
  
  \medskip
  
  Since $G$ is balanced, we get $|$Card$(I) - $Card$(J)|\leq 1$ so condition (i) is satisfied.
  
  \medskip
  
  For every $i\in I$, we denote $$p_i = \mathrm{Card} \{k\in I : v_iv_k \in E\},  q_i = \mathrm{Card} \{k\in J : v_iv_k\in E\},$$ and for every $j\in J$, we denote  
 $$m_j = \mathrm{Card} \{k\in I : v_jv_k\in E\},  n_j = \mathrm{Card} \{k\in J : v_jv_k\in E\}$$  
  
  \medskip
  
  Then $\displaystyle \sum _{i\in I} q_i = \displaystyle \sum _{j\in J} m_j = $Card$(E\backslash (E_1\cup E_2))$. On the other hand, since $G$ is balanced, we have $\displaystyle \sum _{i\in I} p_i = 2 $Card$(E_1), \displaystyle \sum _{j\in j} n_j = 2 $Card$(E_2)$.
  
  \medskip
  
  Then $|\displaystyle \sum _{k\in I}d_k - \displaystyle \sum _{k\in J}d_k| = |\displaystyle \sum _{k\in I}(p_k + q_k) - \displaystyle \sum _{k\in J}(m_k + n_k)| = 2|$Card$(E_1)-$Card$(E_2)|\leq 2$. Thus condition (ii) is also satisfied.
  
  \medskip
  
  To prove the converse, assume conditions (i) and (ii) are satisfied. We define the coloring $c:V\rightarrow \{1,2\}$ as follows: for every $i\in I$ we set $c(v_i) = 0$ and for every $j\in J$ we set $c(v_j) = 1$.
  
 \medskip
 
 Then we have \ Card$(E_1) = \frac{1}{2} \displaystyle \sum _{i\in I} p_i, \ $Card$(E_2) = \frac{1}{2} \displaystyle \sum _{j\in J} n_j, \ \mathrm{and} \\ \displaystyle \sum _{i\in I} q_i = \displaystyle \sum _{j\in J} m_j = $Card$(E\backslash (E_1\cup E_2))$. 

\medskip

On the other hand, $$\displaystyle \sum _{k\in I}d_k = \sum _{k\in I} (p_i+q_i) \ \mathrm{and} \ \displaystyle \sum _{k\in J}d_k = \sum _{k\in J}(m_j+n_j)$$ 
 
 \medskip
 
  Then by condition (ii), we get $|$Card$(E_1) - $Card$(E_2)| = \frac{1}{2}|\displaystyle \sum _{k\in I}d_k - \displaystyle \sum _{k\in J}d_k| \leq 1$.  $\blacksquare $

  \bigskip   
  
  \begin{cor} It is proved in [LLT] that an $r$-regular finite simple graph with $n$ vertices is balanced iff $n$ is even or $r=2$. This fact also follows immediately from Lemma \ref{lem:trivial}. In [KLST], the authors deduce the same fact from their characterization of balanced graphs.
  \end{cor}
  
  \medskip
  
  Lemma \ref{lem:trivial} shows that the balancedness of a graph totally depends on the degree sequence of it. This is no longer the case for $k$-balanced graphs for $k\geq 3$. In fact, the trees $G_1$ and $G_2$ in Figure 3 have the same degree sequence $(1,1,1,1,1,1,1,1,1,1,1,2,2,2,11)$, and it is not difficult to see that $G_1$ is 3-balanced while $G_2$ is not. 
 	
 	\medskip
 	
 	\begin{figure}[h!]
  \includegraphics[width=3in,height=3in]{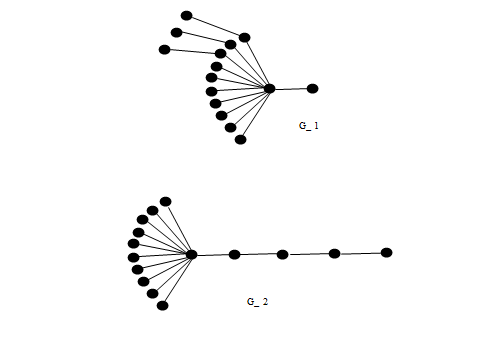}
\caption{The trees $G_1$ and $G_2$ have the same degree sequence; $G_1$ is 3-balanced while $G_2$ is not.}
\label{labelname}
\end{figure}
 	
 	The fact that, for $k\geq 3$, the $k$-balancedness is not determined by the degree sequence causes difficulties in proving that random graphs are $k$-balanced. It also seems plausible that, generically, $k$-balancedness is a weaker condition than balancedness, although it does not seem easy to describe (with a good sufficient condition) when exactly is this true. It is useful to point out the following simple fact.
 	
 	\medskip
 	
 	\begin{prop} For all distinct $m, n\geq 2$ there exists a finite simple graph which is $m$-balanced but not $n$-balanced.
 	\end{prop} 
	
	{\bf Proof.} Let $p$ be a prime number such that $p > \max\{m,n\}$.
	
	\medskip
	
	Let us first assume that $m > n$. If $n$ divides $m$, then the graph $K_{m+1}$ is $m$-balanced but not $n$-balanced. If $n$ does not divide $m$ then the graph $K_{mp}$ is $m$-balanced but not $n$-balanced.
	
  \medskip
  
  Now assume that $m < n$. Then the graph $K_{mp}$ is $m$-balanced but not $n$-balanced. $\blacksquare $

 \bigskip
     
 \section {Combinatorial Lemmas}

 Let $M_n = \{\bar{d} = (d_1, \ldots , d_n) : d_i\in \mathbb{N}, 1\leq d_i\leq n, 1\leq i\leq n,\}$. The elements of $M_n$ consist of sequences of positive integers of length $n$ such that no term is bigger than $n$. We denote $\max (\bar{d}) = \displaystyle \max _{1\leq i\leq n}d_i$. 
 
 \medskip

 Now we introduce the notion of balanced sequences
 
 \medskip
 
 {\bf Definition 3.} A sequence (an element) $\bar{d}\in M_n$ is called {\em balanced} if and only if there exists a partition $\{1,\ldots , n\} = I\sqcup J$ such that
  
  (a) $|$Card$(I)-$Card$(J)|\leq 1$
  
  (b) $|\displaystyle \sum _{k\in I}d_k - \displaystyle \sum _{k\in J}d_k|\leq 2$
	
	\medskip
	
	 The partition $\{1,\ldots , n\} = I\sqcup J$ will be called {\em a balanced partition.}
   
 \medskip
 
  In these new terms, Lemma \ref{lem:trivial} states that a graph is balanced if and only if its degree sequence is balanced. 
  
	\medskip
	
  When the sequence is not balanced, we would like to measure how far it is
from being balanced.

 \medskip
 
  \begin{defn}  Let   $\bar{a} = (a_1, \ldots , a_n)$ be any finite
sequence of non-negative integers. The quantity $$F(\bar{a}) = \displaystyle \min
_{\{1,\ldots ,n\} = I\sqcup J, |\mathrm{Card} (I)-\mathrm{Card}(J)|\leq 1} \left |\displaystyle \sum _{k\in I}a_k -
\displaystyle \sum_{k\in J}a_k\right |$$ will be called {\em the balance of} $\bar{a}$. 
\end{defn}

\medskip

  \begin{rem} By Lemma \ref{lem:trivial}, a sequence $\bar{d}\in M_n$ is balanced if and only if $0\leq F(\bar{d}) \leq 2$. The quantity $F(\bar{a})$, somewhat roughly, measures how far the sequence is from being balanced. For an example, let $n = 8$ and $\bar{a} = (1, 3, 12, 2, 1, 1, 4, 3)$ be a sequence of length 8. It is easy to see $F(\bar{a}) = |(12+1+1+1)-(2+3+3+4)| = 3$.  
  \end{rem}
  
  \medskip

  The following easy lemma will be useful

\medskip

   \begin{lem}\label{lem:lemma1} Let $\bar{a} = (a_1, \ldots , a_n)$ be any finite sequence
of non-negative integers. Then  $F(\bar{a}) \leq \mathrm{max}(\bar{a})$.
\end{lem}

\medskip

{\bf Proof.} We will present a constructive proof.

\medskip

  Without loss of generality, we may assume that $d_1\leq
d_2\leq \ldots \leq d_n$. First, let us assume that $n$ is even, so let $n
= 2m$. We will build two subsets $I, J$ of $\{1,\ldots, n\}$ inductively
such that $\{1,\ldots ,n\} = I\sqcup J, |I| = |J|,$ and $ |\sum _{k\in I}a_k -
\sum_{k\in J}a_k|\leq \mathrm{max}(\bar{a})$.

\medskip

 We divide the sequence into pairs $(d_1,d_2), \ldots ,
(d_{2m-1},d_{2m})$, and we will abide by the rule that exactly one element
of each pair belongs to $I$ and the other element belongs to $J$. We
start by letting $I_1 = \{d_{2m}\}, J_1 = \{d_{2m-1}\}$. Assume now  we
have built the subsets $I_k, J_k, 1\leq k\leq m-1$ such that  $\{d_{2m},
d_{2m-1}, \ldots , d_{2m-2k+2}, d_{2m-2k+1}\} = I_k\sqcup J_k$ and
$|\{d_{2m-2i-2}, d_{2m-2i-1}\}\cap I_k| = 1$ for all $1\leq i\leq k$.

\medskip

Let $S(I_k) = \sum _{i\in I_m}a_i,    S(J_k) = \sum _{I\in J_m}a_i$. If
$S(I_k) > S(J_k)$ then we let $I_{k+1} = I_k\sqcup \{d_{2m-2k-1}\},
J_{k+1} = J_k\sqcup \{d_{2m-2k}\}$ but if  $S(I_k) \leq S(J_k)$ then we
let $I_{k+1} = I_k\sqcup \{d_{2m-2k}\}, J_{k+1} = J_k\sqcup
\{d_{2m-2k-1}\}$, and we proceed by induction. Then we let $I = I_m, J =
J_m$. Clearly, we have $F(\bar{a}) \leq |\sum _{k\in I}a_k -
\sum_{k\in J}a_k| \leq \mathrm{max}(\bar{a})$.

\medskip

If $n$ is odd, then we may replace $\bar{a}$  by $\bar{a'} = (0, a_1,
\ldots , a_n)$ and apply the previous argument. $\blacksquare $
  
\bigskip

 We will need the following notations

\medskip

 \begin{defn} Let $\bar{d} = (d_1, \ldots ,d_n)\in M_n$. We will denote 
$$u(\bar{d}) = \{1\leq i \leq n \ | \ d_i = 1\}, \ v(\bar{d}) = \{1\leq i
\leq n \ | \ d_i = 2\}$$ 
\end{defn}

 \medskip
 
 \begin{lem}\label{lem:lemma2} Let $\bar{d} = (d_1, \ldots , d_n)\in M_n$ such that $|u(\bar{d})|\geq \mathrm{max}(\bar{d})$ and $|v(\bar{d})|\geq \mathrm{max}(\bar{d})$. Then $\bar{d}$ is balanced.
 \end{lem}
  	
\bigskip

{\bf Proof.} Let $\mathrm{max}(\bar{d}) = m$. Without loss of generality we may assume that $d_1 = \ldots = d_m = 1, d_{m+1} = \ldots = d_{2m} = 2$. If $n=2m$ then $\bar{d}$ is clearly balanced so let $n > 2m$ and let $\bar{d}' = (d_{2m+1}, \ldots , d_n)$.

 \medskip

 By Lemma \ref{lem:lemma1}, $F(\bar{d}') \leq m$ hence there exists a partition $\{d_{2m+1}, \ldots , d_n\} = I'\sqcup J'$ such that $|$Card$(I')-$Card$(J')|\leq 1$ and $|\sum _{k\in I'}d_k - \sum _{k\in J'}d_k|\leq m$. Then there exists a partition $\{d_{1}, \ldots , d_{2m}\} = I''\sqcup J''$ such that Card$(I'') =$Card$(J'')$ and $|(\sum _{k\in I''}d_k - \sum _{k\in J''}d_k) - (\sum _{k\in I'}d_k - \sum _{k\in J'}d_k)|\leq 2$. By letting $I = I'\sqcup I'', J = J'\sqcup J''$ we obtain that $\{1, \ldots , n\} = I\sqcup J$, $|$Card$(I)-$Card$(J)|\leq 1$, and $|\sum _{k\in I}d_k - \sum _{k\in J}d_k|\leq 2$. $\blacksquare $

 \bigskip

 \section{Proof of Theorem A}

	First, we will discuss the case of labeled trees. The following theorem of J.W.Moon will play a crucial role   

  \medskip
  
  \begin{thm}[See (M)]\label{thm:moon} If $\epsilon > 0$ is a fixed positive constant, then in a random labeled tree with $n$ vertices, the maximal degree $d_{\max }$ satisfies the following inequality $$(1-\epsilon)\frac{\log n}{\log \log n} < d_{\max } < (1+\epsilon)\frac{\log n}{\log \log n}$$
  \end{thm}
  
  \medskip
  
  \begin{rem}\label{rem:bound} By choosing $\epsilon = 0.1$ we obtain that 
   $$0.9\frac{\log n}{\log \log n} < d_{\max } < 1.1\frac{\log n}{\log \log n}$$ in a random tree with $n$ vertices. 
   \end{rem}
     
 \bigskip
 
  We will use only the upper bound in the inequality of Remark \ref{rem:bound}. Besides the upper bound on the maximal degree in random trees, we also need a lower bound on the number of vertices with degree 1, and with degree 2. Notice that, since the sum of degrees of a tree with $n$ vertices is exactly $2n-2$, at least half of the vertices have degree either 1 or 2. However, we need a linear lower bound for the number of vertices of degree 1 and for the number of vertices of degree 2 separately.
		
	\medskip
		
 Let $X_i(T), 1\leq i\leq 2$ be the random variable which denotes the number of vertices of degree $i$ in a labeled tree $T$ with $n$ vertices. Also let $\mu = \frac {n}{e}, \sigma _1^2 = \frac {n}{e}(1- \frac {2}{e}), \sigma _2^2 = \frac {n}{e}(1- \frac {1}{e})$. It has been proved by A.R\'{e}nyi (see [R]) that the asymptotic distribution of random variable $\frac {X_1-\mu }{\sigma _1}$ is normal with mean $\mu $ and variance $\sigma _1^2$. A similar result has been proved for the random variable $\frac {X_2-\mu }{\sigma _2}$, by A.Meir and J.W.Moon (see [MM]), namely, that the asymptotic distribution of the random variable $\frac {X_2-\mu }{\sigma _2}$ is normal with mean $\mu $ and variance $\sigma _2^2$. Combining these two results we can state the following theorem (due to A.R\'{e}nyi and A.Meir-J.W.Moon)

 \medskip

\begin{thm}\label{thm:meir} Let $\alpha , \beta $ be fixed real numbers, $\alpha < \beta $; and for $i\in \{1,2\}$, let $P_i(\alpha , \beta )$ denotes the probability that $\alpha < \frac {X_i-\mu }{\sigma _1} < \beta $. Then $$\lim _{n\rightarrow \infty }P_i(\alpha , \beta ) = \frac {1}{\sqrt{2\pi }}\int _{\alpha }^{\beta }e^{-\frac{1}{2}t^2}dt$$
\end{thm}

\bigskip

 We need the following immediate corollary of this theorem
 
 \medskip

 \begin{cor}\label{cor:inequality} In a random labeled tree with $n$ vertices, for all $i\in \{1,2\}$, $X_i\geq 2\frac{\log n}{\log \log n}$.
\end{cor}

\medskip

  Now, in the case of random labeled trees, Theorem A immediately follows from Theorem \ref{thm:moon}, Corollary \ref{cor:inequality}, and Lemma \ref{lem:lemma2}.

 \bigskip

  {\bf {\em The case of unlabeled trees:}} We will use the results analogous to Theorem \ref{thm:moon} and Theorem \ref{thm:meir}. The analogue of Theorem \ref{thm:moon} is proved by W.Goh and E.Schmutz:
  
  \medskip
  
  \begin{thm}[See (GS)]\label{thm:bounds} There exists positive constants $c_1, c_2$ such that in a random unlabeled tree with $n$ vertices the maximum degree $d_{\mathrm{max}}$ satisfies the inequality $c_1\mathrm{log}(n) < d_{\mathrm{max}} < c_2\mathrm{log}(n)$.
  \end{thm}
  
  \medskip

  Now, for any $k\in \mathbb{N}$ let the random variable $Y_k$ denotes the number of vertices of degree $k$ in a random unlabeled tree with $n$ vertices. The following theorem is due to M.Drmota and B.Gittenberger; in the case of $k\in \{1,2\}$, as a special case, it provides an analogue of Theorem \ref{thm:meir}.
		
		\medskip
		
	\begin{thm}[See (DG)] For arbitrary fixed natural $k$, there exists positive constants $\mu _k $ and $\sigma _k$ such that the limiting distribution of $Y_k$ is normal with mean $\mu (n) \sim \mu _kn$ and variance $\sigma (n) \sim \sigma _k^2n$.
	\end{thm}
			
	\medskip
	
	\begin{cor}\label{cor:linearbound}For all $c>0$ and $i\in {1,2}$, in a random unlabeled tree with $n$ vertices $Y_i > c\mathrm{log}(n)$.	
	\end{cor}
   
\bigskip

 Now, in the case of unlabeled trees, the claim of Theorem A follows from Theorem \ref{thm:bounds}, Lemma \ref{lem:lemma2}, and Corollary \ref{cor:linearbound}.
 
 \bigskip
  
  \section{$k$-balanced trees: proof of Theorem B}

  In this section we will assume that $k\geq 3$. The fact that the $k$-balancedness is not determined by the degree sequence causes significant difficulties in proving that random graphs are balanced. We nevertheless prove that random trees are strongly $k$-balanced by more careful study of $k$-balancedness.

	\medskip
	 
	First, we need to prove the following technical lemma.
	
	\medskip
	
	\begin{lem}\label{lem:w=2} Let $G = (V,E)$ be a tree and $u,v$ be distinct vertices of $G$ with degrees at least $\frac{|G|}{3}$. Let also  $p,q$ be distinct pre-leaf vertices of $G$. Then there exists a strongly 3-balanced coloring $c:V\rightarrow \{1,2,3\}$ of $G$ such that $c(u)\neq c(v)$ and $c(p)\neq c(q)$.  
	\end{lem}
	
	\medskip
	
	{\bf Proof.} The proof is by induction on $n = |G|$. For $n\leq 5$ the claim is obvious so we will assume that $n\geq 6$ and the claim holds for all trees of order less than $n$.
	
	\medskip
	
	Assume that at least one of the following two conditions hold:
	
	(c1) there exists $z\in \{p,q\}\backslash \{u,v\}$ such that $deg(z)\geq 3$;
	
	(c2) there exists a leaf vertex not adjacent to any of the vertices $u, v, p, q$.
	
	Then there exists a leaf $w$ such that if $G'$ is a complete subgraph on $V\backslash \{w\}$, then, in the tree $G'$, we have $\min \{deg(u), deg(v)\} \geq \frac{|G'|}{3}$, and $p, q$ are still pre-leaf vertices.
	
	\medskip
	
	By inductive hypothesis, there exists a strongly 3-balanced coloring $c_0:V\backslash \{w\}\rightarrow \{1,2,3\}$ of $G'$ such that $c_0(u)\neq c_0(v)$ and $c_0(p)\neq c_0(q)$. Let $w_0$ be the unique vertex of $G$ adjacent to $w$. Without loss of generality, we may assume that $c_0(w_0) = 1$ and $|c_0^{-1}(2)|\leq |c_0^{-1}(3)|$.  
	
	\medskip
	
	If $|c_0^{-1}(1)|\geq |c_0^{-1}(2)|$ then we let $c(w) = 2$ thus extending $c_0$ to a strongly 3-balanced coloring $c:V\rightarrow \{1,2,3\}$ of $G'$ such that $c(u)\neq c(v)$ and $c(p)\neq c(q)$.
	
	\medskip
	
	If, however, $|c_0^{-1}(1)| < |c_0^{-1}(2)|$ then there exists $r\in \{u,v\}$ such that $c_0(r)\neq 1$; also, since $deg(r)\geq \frac {|G|}{3}$, there exists a branch $B$ of $G'$ with respect to $r$ which is disjoint from $c_0^{-1}(1)$. Let $x$ be a leaf vertex in $B$. Then $x\notin \{u,v,p,q\}$ and $c_0(x) \neq 1$. We define $c:V\rightarrow \{1,2,3\}$ as follows:
	
	\begin{displaymath}  c(z) =  \left\{\begin{array}{lcr} c_0(z) &  \mathrm{if} \ z\in V\backslash \{w, x\} \\ 1 &  \mathrm{if} \ z = x \\  c_0(x) &  \mathrm{if} \ z = w  \end{array}   \right. \end{displaymath}
	
	Notice that because of the inequality  $|c_0^{-1}(1)| < |c_0^{-1}(2)|\leq |c_0^{-1}(3)|$,  we have $|c_0^{-1}(2)| = |c_0^{-1}(3)|$ and $|c_0^{-1}(1)| = |c_0^{-1}(2)|-1$. Then the map  $c:V\rightarrow \{1,2,3\}$ is a strongly 3-balanced coloring.
	
	\medskip
	
	Now, suppose that none of the conditions (c1) and (c2) hold. Let $P$ be the path in $G$ starting at $u$ and ending at $v$ (it may possibly consist of just the vertices $u$ and $v$). Then the tree $G$ satisfies the following conditions: there exists two vertices $z_1, z_2$ in $P$ and paths $R_1, R_2$ starting at $z_1, z_2$ respectively such that any vertex of $G$ either belongs to one of the paths $P, R_1, R_2$ or it is a leaf vertex adjacent to one of the vertices $u, v$. Then it is straightforward to build a strongly 3-balanced coloring $c:V\rightarrow \{1,2,3\}$ satisfying the conditions $c(u)\neq c(v)$ and $c(p)\neq c(q)$.  $\blacksquare $
	   
	\medskip
	
	The following proposition is interesting in itself; it will also play a key role in proving Theorem B.
	
	\medskip
	
	\begin{prop}\label{prop:k=3} Let $G=(V,E)$ be a tree with $n$ vertices where $d_{max}(G) \leq \frac{n}{3}$. Then $G$ is strongly $3$-balanced. Moreover, for any two distinct pre-leaf vertices $p$ and $q$ of $G$ there exists a strongly 3-balanced coloring $c:V\rightarrow \{1,2,3\}$ such that $c(p)\neq c(q)$. 
	 \end{prop}
	 
	 \medskip
	 
	 {\bf Proof.} The proof will be by induction on $n$. For $n\geq 8$ we have $d_{max}(G)\leq 2$ hence $G$ is isomorphic to a string, thus the claim is obvious. Let us now assume that $n\geq 9$, and the claim holds for all trees $G'$ of order less than $n$ with $d_{max}(G')\leq \frac{|G'|}{3}$.
	 
	 \medskip
	 
	  Let $G =(V,E)$ and $n = 3k+r, r\in \{0,1,2\}$. We will consider the following three cases separately:
	  
	  \medskip
	  
	  {\em Case 1.} $r = 1$.
	  
	  Let $v$ be a leaf of $G$, $V' = V\backslash \{v\}$, and let $G' = (V', E')$ be the complete subgraph of $G$ on $V'$. Then we have $$d_{max}(G')\leq d_{max}(G) \leq k \leq \frac{|G'|}{3}.$$ By inductive hypothesis, there exists a strongly 3-balanced coloring $c':V'\rightarrow \{1,2,3\}$ of $G'$. 
	  
	  \medskip
	  
	  On the other hand, $v$ is adjacent to exactly one vertex in $G$; let $u$ be this vertex. Let $j$ be any element of $\{1,2,3\}\backslash \{c'(u)\}$. We extend the coloring $c'$ of $G'$ to a strongly 3-balanced coloring $c:V\rightarrow \{1,2,3\}$ by defining $c(v) = j$. 
	  
	  \medskip
	  
	  {\em Case 2.} $r=2$.
	  
	  Let $v_1, v_2$ be distinct leaves and $u_1, u_2$ be the only vertices of $G$ adjacent to $v_1, v_2$ respectively ($u_1$ and $u_2$ are not necessarily distinct). Let also $G'$ be the complete subgraph of $G$ on the set $V\backslash \{v_1, v_2\}$. Then we still have the inequality $d_{max}(G')\leq d_{max}(G) \leq k \leq \frac{|G'|}{3}$. Hence, by inductive assumption, there exists a strongly 3-balanced coloring $c':V'\rightarrow \{1,2,3\}$ of $G'$.
	  
	  \medskip
	  
	  Then there exist distinct $j_1, j_2\in \{1,2,3\}$ such that $j_1\neq c'(u_1)$ and $j_2\neq c'(u_2)$. Thus we can extend $c'$ to a strongly 3-balanced coloring of $G$ by defining $c(v_1) = j_1$ and $c(v_2) = j_2$. 
	  
	  \medskip
	  
	  {\em Case 3.} $r=0$. 
	  
	  The major difference in this case compared with the previous two cases is that when we obtain $G'$ by deleting some arbitrary three leaves $v_1, v_2, v_3$ from $G$ ($G$ possesses three leaf vertices unless it is isomorphic to a string) we may loose the inequality $d_{max}(G')\leq \frac{|G'|}{3}$. Also, suppose $u_1, u_2, u_3$ are the vertices adjacent to $v_1, v_2, v_3$ respectively ($u_1, u_2, u_3$ are not necessarily distinct).  If we have the inequality $d_{max}(G')\leq \frac{|G'|}{3}$ then by inductive assumption we would have a strongly 3-balanced coloring $c':V\backslash \{v_1, v_2, v_3\}\rightarrow \{1,2,3\}$, however, if $c'(u_1) = c'(u_2) = c'(u_3)$ then it becomes problematic to extend $c'$ to a strongly 3-balanced coloring $c:V\rightarrow \{1,2,3\}$. Thus we need to employ different and more careful tactics. 
	  
	  \medskip
	  
	  We will prove the following lemma which suffices for the proof of Proposition \ref{prop:k=3} in the case $r=0$.
	  
	  \begin{lem} Let $G = (V, E)$ be a tree with $n=3k$ vertices where $d_{max}(G) \leq k$, and let $p, q$ be distinct pre-leaf vertices of $G$. Then there exists a strongly 3-balanced coloring $c:V\rightarrow \{1,2,3\}$ such that $c(p)\neq c(q)$.
	  \end{lem}
	  
	  {\bf Proof.} The proof of the lemma will be again by induction on $k$. The ``$c(p)\neq c(q)$ part" of the claim will be needed to make the step of the induction. For $k\leq 2$, the graph $G$ is isomorphic to a string thus the claim is obvious. For $k=3$ it can be seen by a direct checking (we leave this to a reader as a simple exercise). So let us assume that $k\geq 4$. 
	  
	   \medskip  
	  
	  Let $W = \{v\in V \ | \ d(v)=k\}$. Let also $deg(p)\leq deg(q)$. We will consider the following cases (the notations in each case will be independent of the notations of other cases):
	  
	  \medskip
	  
	  {\bf Case A:} $W = \emptyset $ and $p$ is not special.
	  
	  \medskip
	  
  Let $v_1, v_2, v_3$ be distinct leaves such that $v_1$ is adjacent to $p$, $v_2$ is adjacent to $q$, and $v_3$ is adjacent to a vertex $w$ distinct from  $p$ and $q$. We let $G'$ be the complete subgraph on $V\backslash \{v_1, v_2, v_3\}$. Then $|G'| = 3(k-1)$ and we have $d_{max}(G')\leq k-1$. By inductive hypothesis, there exists a strongly 3-balanced coloring $c_0:V\rightarrow \{1,2,3\}$ such that $c_0(p)\neq c_0(q)$. Without loss of generality we may assume that $c_0(p) = 1, c_0(q) = 2$. Then we extend $c_0$ to a strongly 3-balanced coloring $c:V\rightarrow \{1,2,3\}$ as follows: if $c_0(w)\in \{1,2\}$ then we let $c(v_1) = 2, c(v_2) = 1, c(v_3) = 3$; and if $c_0(w) = 3$ then we let $c(v_1) = 3, c(v_2) = 1, c(v_3) = 2$.       	  
	  
	  \bigskip
	  
	  {\bf Case B:} $W = \emptyset $ and $p$ is special.
	  
	  \medskip
	  
  Let $v_1$ be the only leaf adjacent to $p$, $u$ be the unique non-leaf vertex adjacent to $p$, $v_2$ be a leaf vertex not adjacent to $u$, and $w$ be the unique vertex adjacent to $v_2$. We let $G'$ be the complete subgraph on $V\backslash \{v_1, v_2, p\}$. Then $|G'| = 3(k-1)$ and $d_{max}(G')\leq k-1$. By inductive hypothesis, there exists a strongly 3-balanced coloring $c_0:V\rightarrow \{1,2,3\}$. Then we extend $c_0$ to a strongly 3-balanced coloring $c:V\rightarrow \{1,2,3\}$ as follows:  we let $c(p)\in \{1,2,3\}$ such that $c(p)$ is distinct from $c_0(u)$ and $c_0(q)$. Then we define $c(v_2)\in \{1,2,3\}$ such that $c(v_2)$ is distinct from $c_0(w)$ and $c(p)$. Finally we let $c(v_1)\in \{1,2,3\}$ such that $c(v_1)$ is distinct from $c(p)$ and $c(v_2)$. Notice also that we obtain $c(p)\neq c(q)$.

	  \bigskip
	   
	  {\bf Case C:} $|W| = 1, W = \{v_0\}, deg(p)\geq 3$ and there exists a leaf vertex adjacent to $v_0$.
	  
	  \medskip
	  
	 This case is similar to  Case A. Since $|W| = 1$ and $deg(p)\leq deg(q)$, we have $p\neq v_0$. If $q\neq v_0$, we let $v_1, v_2, v_3$ be leaves adjacent to $p, q, v_0$ respectively; and if $q=v_0$, we let $v_1, v_2$ be leaves adjacent to $p, q$ respectively, and $v_3$ be a leaf not adjacent to any of the vertices $p, q$. We define $G'$ to be the complete subgraph on $V\backslash \{v_1,v_2,v_3\}$. Then $d_{max}(G')\leq \frac{|G'|}{3}$ hence $G'$ admits a strongly 3-balanced coloring $c':V\backslash \{v_1,v_2,v_3\}\rightarrow \{1,2,3\}$ such that $c'(p) \neq c'(q)$. We extend $c'$ to a strongly 3-balanced coloring to $c:V\rightarrow \{1,2,3\}$ as in Case A.
	 
	 \bigskip
	 
	 {\bf Case D:} $|W| = 1, W = \{v_0\}$, $p$ is special and there exists a leaf vertex adjacent to $v_0$.
	 
	 \medskip
	 
	 This case is similar to  Case B. Let $v_1$ be the only leaf adjacent to $p$, $u$ be the unique non-leaf vertex adjacent to $p$, $v_2$ be a leaf vertex adjacent to $v_0$. We let $G'$ be the complete subgraph on $V\backslash \{v_1, v_2, p\}$. Then $|G'| = 3(k-1)$ and $d_{max}(G')\leq k-1$. By inductive hypothesis, there exists a strongly 3-balanced coloring $c_0:V\backslash \{v_1, v_2, p\}\rightarrow \{1,2,3\}$. Then we extend $c_0$ to a strongly 3-balanced coloring $c:V\rightarrow \{1,2,3\}$ as follows:  we let $c(p)\in \{1,2,3\}$ such that $c(p)$ is distinct from  $c_0(u)$ and $c_0(q)$. Then we define $c(v_2)\in \{1,2,3\}$ such that $c(v_2)$ is distinct from  $c_0(v_0)$ and $c(p)$. Finally we let $c(v_1)\in \{1,2,3\}$ such that $c(v_1)$ is distinct from $c(p)$ and $c(v_2)$.

	 \bigskip
	  
	  {\bf Case E:} $|W| = 1, W = \{v_0\}$, and there is no leaf vertex adjacent to $v_0$.
	  
	  \medskip
	  
	  Then, necessarily, there exists a special vertex $v$ adjacent to $v_0$. Let $v_1$ be the unique leaf adjacent to $v$. Let also $v_2$ be a leaf not adjacent to any of the vertices $p,q,v$ (such a leaf exists because $k\geq 4$), and let $w$ be the unique vertex adjacent to $v_2$.
	 
	  \medskip
	  
	  We define $G'$ to be the complete subgraph on $V\backslash \{v, v_1, v_2\}$. By inductive assumption, there exists a strongly 3-balanced coloring $c_0:V\backslash \{v, v_1, v_2\}\rightarrow \{1,2,3\}$, moreover, if $p, q\in V\backslash \{v, v_1, v_2\}$ then $c_0(p)\neq c_0(q)$.

	  \medskip
	  
	  If $\{p, q\}\subset V\backslash \{v, v_1, v_2\}$, then we let $c(v_2)$ be any element of $\{1,2,3\}$ distinct from $c_0(w)$. Then we let $c(v)$ be any element of $\{1,2,3\}$ distinct from $c_0(v_0)$ and $c(v_2)$. Finally, we let $c(v_1)$ be any element of $\{1,2,3\}$ distinct from $c(v)$ and $c(v_2)$. Thus we have extended $c_0$ to a strongly 3-balanced coloring $c:V\rightarrow \{1,2,3\}$ such that $c(p)\neq c(q)$.
	  
	  \medskip
	  
	  If $\{p, q\}\cap \{v, v_1, v_2\} \neq \emptyset $ then $\{p, q\}\cap \{v, v_1, v_2\} = \{v\}$ and we may assume that $p = v$. Then we let $c(v)$ be any element of $\{1,2,3\}$ distinct from $c_0(v_0)$ and $c_0(q)$; then we let $c(v_2)$ be any element of $\{1,2,3\}$ distinct from $c_0(w)$ and $c(v)$; finally we let $c(v_1)$ be any element of $\{1,2,3\}$ distinct from $c(v)$ and $c(v_2)$.

	  \bigskip

	  {\bf Case F:} $|W| \geq 2$. 
	  
	  \medskip
	  
	 In this case the claim follows immediately from Lemma \ref{lem:w=2}. $\blacksquare $
	  
	 \bigskip
	 
	 Now we can prove an analogous result for $k$-balanced graphs.
	 
	 \begin{prop}\label{prop:k-balanced} Let $G = (V,E)$ be a tree with $n$ vertices where $d_{max}(G) \leq \frac{n}{k}$ and $k\geq 3$. Then $G$ is strongly $k$-balanced.
	 \end{prop}
	 
	 {\bf Proof.} The proof is by induction on $k$. For $k=3$, the claim is true by Proposition \ref{prop:k=3}. 
	 
	 \medskip
	 
	 Assume now $k\geq 4$. Then the tree $G$ has $m = \lfloor \frac{n}{k} \rfloor $ vertices $v_1, \ldots , v_m$ such that $d(v_i)\leq 2, 1\leq i\leq m$, moreover, for all distinct $i,j\in \{1, \ldots , m\}$, the vertices $v_i$ and $v_j$ are not connected by an edge. Let also $V_0 = \{v_1, \ldots , v_m\}$, and $G_1$ be a complete graph on the subset $V\backslash V_0$. Then $G_1$ is a forest with $n-m$ vertices but with $d_{max}(G_1) \leq d_{max}(G)$. Then $G_1$ is a subgraph of a tree $G_2$ with $n-m$ vertices where $d_{max}(G_2)\leq d_{max}(G)$.
	 
	 \medskip
	 
	Then $d_{max}(G_2) \leq d_{max}(G)\leq \frac{n}{k} = \frac{1}{k-1}(n-\frac {n}{k})\leq \frac{1}{k-1}(n-m) \leq \frac{|G_2|}{k-1}$. Then, by inductive hypothesis, $G_2$ is strongly $(k-1)$-balanced, hence $G_1$ is strongly $(k-1)$-balanced. Since no two elements of $V_0$ are adjacent, we obtain that $G$ is strongly $k$-balanced. $\blacksquare $  
	
	\medskip
	
	Now, for random labeled trees, Theorem B follows immediately from Theorem \ref{thm:moon} and Proposition \ref{prop:k-balanced}; and for random unlabeled trees, it follows immediately from Theorem \ref{thm:bounds} and Proposition \ref{prop:k-balanced}.

  \vspace{1cm}
  
  {\bf R e f e r e n c e s:}
  
  \bigskip
	
	[BK] Ben-Eliezer, I. and Krivelevich, M. \ Perfectly balanced partitions of smoothed graphs. \   
 {\em Electronic Journal of Combinatorics}, {\bf 16} (1) (2009), note N14.   
	
	\medskip
	
	[B] Bollob\'{a}s, B. \ Random Graphs. Academic Press Inc. (1985)
	
	\medskip
	
	[BG] Bollob\'{a}s, B. and Guy, R. \ Equitable and proportional coloring of trees, \ {\em Journal of Combinatorial Theory, Series B} {\bf 34}, (2): (1983) 177–186
	
	\medskip
	
	[Cah1] Cahit, I. \ Cordial graphs: a weaker version of graceful and harmonious graphs, \ {\em Ars Combinatoria}, {\bf 23} (1987) 201-207.
	
	\medskip
	
	[Cah2] Cahit, I. \ On cordial and 3-equitable graphs, \ {\em Utilitas Mathematica}, {\bf 37} (1990) 189-198.
	
	\medskip
	
  [Cay] Cayley, A. \ A Theorem on Trees, {\em Quart.J.Math}, {\bf 23}, (1889) 376-378.
  
  \medskip
  
  [DG] Drmota, M.  Gittenberger, B. \  Distribution of nodes of given degree in random trees.  \ {\em Journal of Graph Theory,} \ {\bf 31/3} (1999), 227-253.  
	
	\medskip
	
	[Ga] Gallian, J.A. \ A dynamical survey of graph labeling. \ {\em The Electronics Journal of Combinatorics}, {\bf 16} (2009)
	
	\medskip
	
	[GS] Goh, W. Schmutz, E. \ Unlabeled trees: distribution of the maximum degree. \ {\em Random Structures and Algorithms}, {\bf 5(1)}, (1994), 13-24.
	
	\medskip
	
	[GS] Graham, R; Sloane, N. \ On additive bases and harmonious graphs. \ {\em SIAM Journal of Algebraic and Discrete Mathematics}, {\bf 1}, (1980) 382-404.
	
	\medskip
  
	[KLST] Kong, M.C., Lee, S-M., Seah, E., Tang, A. \  A  complete characterization of balanced graphs.  (English summary) {\em J. Combin.  Math.  Combin.  Comput.}  {\bf 66} (2008), 225-236.  

	\medskip
		
  [LLT] Lee, S-M., Liu, A., and Tan, S.K. \ On balanced graphs, {\em Congr. Numer.}, {\bf 87} (1992), 59-64.
   
  \medskip
  
  [M] Moon, J.W. On the maximum degree in a random tree. \ {\em Michigan Math. Journal}, Volume 15, Issue 4, (1968), 429-432. 
	
	\medskip
	
	[MM] Meir, A., Moon, J.W. \ On nodes of degree two in random trees. \ {\em Mathematika}, iss.2, {\bf vol. 15} (1968), 188-192. 
	
	\medskip
	
	[O] Otter, R. \ The number of trees. \ {\em Annals of Mathematics}, Second Series, {\bf 49} (3), 583-599.
	
	\medskip
	
	[W] West, D. B., \ {\em Introduction to Graph Theory}, 2nd ed.  
   Prentice-Hall, Inc., Upper Saddle River NJ (2001), 82-83.

 \end{document}